\documentclass[10pt,a4paper]{article}
\usepackage[cp866]{inputenc}
\usepackage{amsmath,amsfonts,amscd,amssymb,latexsym}

\def\rank{\mathop{\fam 0 rank}\nolimits}
\def\sgn{\mathop{\fam 0 sgn}\nolimits}
\def\codim{\mathop{\fam 0 codim}\nolimits}
\def\Ann{\mathop{\fam 0 Ann}\nolimits}
\newcommand{\dve}[1]{\stackrel{\wedge}{#1}}

\begin{document}


\begin{center}
{\Large\bf
Symmetric power of the Grassmann variety}
\end{center}

\begin{center}
{\large V.Yu.Gubarev}
\end{center}

{\it Abstract.}

In this work, we study the subspace $V_0$ of $S^{m}(\wedge^{k}\mathbb{R}^{n})$
spanned by elements $(x_1\wedge\ldots \wedge x_k)^m$, where
$x_i\in\mathbb{R}^{n}$. The problem of finding its
base and dimension is solved. The same problem, but expressed in terms
of polynomials from matrix minors, was initially solved by W. Hodge.
The new result of this paper is the explicit base of $V_0$ and considerable
simplification of the formula for $\dim V_0$.

We also find a set of linear and quadratic relations of rank 3 and 4
defining the Grassmann variety $G_{n,k}$ in the space
$S^{2}(\wedge^{k}\mathbb{R}^{n})$.

\medskip


\begin{center}
\bf \large Introduction
\end{center}

Let $\wedge^{k}\mathbb{R}^{n}$ be the $k$-th exterior power of
$\mathbb{R}^{n}$,
$V(m,n,k)=S^{m}(\wedge^{k}\mathbb{R}^{n})$ be its $m$-th
symmetric power, and let us denote by
$V_0$ the subspace
$L((x_1\wedge\ldots\wedge x_k)^m: x_i\in\mathbb{R}^{n})$.
In this work we study the question
about finding the dimension (denote it by $N(m,n,k)$) and
a base of $V_{0}$.

The space $V(m,n,2)$ appeared in \cite{Shar1} in the definition
of one form of Saint Venant operator. The problem of finding $N(m,n,2)$
came from \cite{Shar2} near the process of recovery of solenoidal part
of symmetric tensor field. The author in \cite{Gub1} solved
this problem in the cases $k=2$ and $n\leq 4$, showing that by
these restrictions $V_0\ne V(m,n,2)$ only in the case $V(m,4,2)$:
Then $N(m,4,2)=\binom{m+5}{5}-\binom{m+3}{5}$. In \cite{Gub2} the problem
was solved for $m=2$ and, in general case, an algorithm of finding base
and dimension was suggested. There were also
some generalizations of the problem to the case of infinite
dimension and the case of anticommutative algebras in \cite{Gub2}.

In $\S1$, we state prerequisite definitions and facts from multilinear
algebra, combinatorics, and algebra of polynomials of
matrix minors. Also, chains of basic vectors are defined there.

In $\S2$, we find base and dimension of $V_0$ in general case.

It is necessary to notice that these results reprove the result
expressed in the terms of polynomials of matrix minors and received
by W.~Hodge in \cite{Hodge} (DeConcini et al. \cite{Ita} stated the result
with by means of Young diagrams technique).

In \cite{Kasman}, a new set of quadratic relations defining projective
presentation of the Grassmann variety and having rank 6 was built.

In $\S3$, we state a set of linear and quadratic relations defining
the Grassmann variety in the space $S^{2}(\wedge^{k}\mathbb{R}^{n})$.
\medskip

\begin{center}
{\bf \large \S1 Main notions }\\
\large 1.1 Decomposable vectors
\end{center}

Let $\wedge^{k}\mathbb{R}^{n}$ be the $k$-th exterior power
of $\mathbb{R}^{n}$, and let $\{e_i,\,i=1,\ldots,\,n\}$ be
some fixed base of $\mathbb{R}^{n}$. We will denote a vector
$e_{i_1}\wedge e_{i_2} \wedge \ldots \wedge e_{i_k}$ as
$e_{\alpha}$, where $\alpha=(i_1,\ldots,i_k)$. Let
$\Pi=\{ \alpha : \alpha=(i_1,\ldots,i_k) : 1\leq i_1<i_2\ldots
<i_k\leq n \}$. It is well known that the set of vectors $\{e_{\alpha} :
\alpha \in \Pi \}$ forms a base of the space $\wedge^{k}\mathbb{R}^{n}$,
and $\dim\,\wedge^{k}\mathbb{R}^{n}=\binom{n}{k}$.

Nonzero vector $v\in \wedge^k\mathbb{R}^n$ is called decomposable if it can be
presented in the form $a_1\wedge\ldots\wedge a_k$, where $a_i\in\mathbb{R}^n$.

Suppose the decomposition of a vector $v\in\wedge^k\mathbb{R}^n$ is known:
\begin{eqnarray}\label{Rasl}
v=\sum_{1\leq i_1<\ldots<i_k\leq n}\mu_{i_1\ldots
i_k}e_{i_1}\wedge \ldots\wedge e_{i_k}.
\end{eqnarray}
Note that if $v=a_1\wedge\ldots\wedge a_k$ and
$A=(a_{ij})$ is a matrix of size $k\times n$  whose rows
come from coordinates of vectors $a_1,\ldots ,a_k$ in the base
$e_1,\ldots,e_n$, then the coefficient $\mu_{i_1\ldots i_k}$
in (\ref{Rasl}) equals the minor of rank $k$
of $A$ formed by columns with numbers $i_1,\ldots,i_k$.

The following theorem states necessary and sufficient conditions
for the vector $v$ to be decomposable:
\medskip

\textbf{Theorem 1 \cite{Vin,McDgl}.}
{\it A vector $v\in
\wedge^{k}\mathbb{R}^n$ is decomposable if and only if the relations
\begin{eqnarray}\label{Plue}
\sum_{p=1}^{k+1}(-1)^{p}\mu_{i_1\ldots \dve{i_{p}}\ldots
i_{k+1}}\mu_{i_{p}j_1\ldots j_{k-1}}=0
\end{eqnarray}
hold for all $1\leq i_1<\ldots<i_{k+1}\leq
n,\,\,\,1\leq j_1<\ldots<j_{k-1}\leq n$. }

\medskip
Relations (\ref{Plue}) are known as {\it Pl\"ucker relations}.
There exists bijective (exact to scalar multiplication) mapping
(see \cite{Vin,Kos})
$$
\Ann(x_1\wedge \ldots \wedge x_k)=L(x_1,\ldots,x_k)
$$
from the set $D$ of decomposable nonzero vectors $\wedge^k\mathbb{R}^n$
into the Grassmann variety $G_{n,k}$ of $k$-dimensional subspaces $\mathbb{R}^n$.
Theorem 1 allows us to present $G_{n,k}$ as algebraic variety in
the projective space $P(\wedge ^{k}\mathbb{R}^n)$.

\begin{center}
\large 1.2 Young tableaux
\end{center}

We are following the book \cite{Fulton} during this paragraph.
Young tableau $\lambda$ with rows $b_1, b_2, ..., b_l$ ($b_1
\geq b_2 \geq \ldots \geq b_l > 0$) is a table whose rows are
ordered by left side and contain correspondingly $b_1, b_2,
\ldots, b_l$ cells. The ordered set of rows lengths generates
a decomposition $(b_1,b_2,\ldots,b_l)$ of the number
$|\lambda|=b_1+\ldots+b_l$.

Each cell of such a table contains a natural number
not exceeding $n$, and these numbers strongly increase inside
every column.

The hook of a cell in a Young tableau is a number of cells
standing below or to the right plus one.
Let us denote the hook of $\alpha$ as $h(\alpha)$.
The content of the cell $\alpha$ is $s(\alpha)=j-i$ if $\alpha$
stands in the $i$-th row and $j$-th column.

If all entries of a Young tableau are distinct, i.e., we have an
arrangement of numbers $1,2,\ldots,|\lambda|$, and
the entries of every its row (strongly) increase then such a Young tableau
is called standard. The well-known hook formula computes a number of standard
Young tableaux
$\lambda$
\begin{eqnarray} \label{HookF}
f^{\lambda} = \frac{|\lambda|}{\prod\limits_{\alpha \in \lambda}
h(\alpha)}.
\end{eqnarray}

If the entries of a Young tableau weakly increase (i.e.,
increase but not strongly) in any row then such a tableau is
called semistandard.

The Stanley formula calculates
the number of semistandard Young table\-aux of the shape $\lambda$,
containing entries not exceeding $n$:
\begin{eqnarray} \label{StanleyF}
d^{\lambda,n} = \frac{\prod\limits_{\alpha \in \lambda}
(n+s(\alpha))}{\prod\limits_{\alpha \in \lambda} h(\alpha)}.
\end{eqnarray}

\begin{center}
\large 1.3 Polynomials in matrix minors
\end{center}

Let $(x_{i j})$ be a matrix of size $k\times n$, whose elements
$x_{ij}$ are independent variables. We denote a minor of rank $s$
by means of its
diagonal elements:
$|x_{i_1 j_1}\ldots x_{i_s j_s}|$ stands for the minor built from
rows $i_1,\dots, i_s$ and columns $j_1,\dots, j_s$.

A $k$-connex (a product of $k$-th degree \cite{HodgePido}) of type
$(l_1,\ldots,l_k)$ is a polynomial in variables $x_{ij}$
which is homogeneous of degree $l_i$
in variables $(x_{i1},\ldots$, $x_{in})$
and can be represented in the form of a polynomial in
the following minors of the matrix $(x_{i j})$:
\begin{equation}\label{Minors}
x_{1 j_1}, |x_{1 j_1}x_{2 j_2}|,\ldots,|x_{1j_1}\ldots x_{k j_k}|,
\quad 1\leq j_1<j_2<\dots<j_k\leq n.
\end{equation}

A $k$-connex for $k=1$ is just a form of degree $l_1$ in
$(x_{1 1},\ldots,x_{1 n})$, and $\binom{n+l_1-1}{l_1}$ equals a number
of linear independent $1$-connexes of type $l_1$. For $k>1$
the different products
$$
(x_{11})^{\rho_{1}} \ldots (x_{1n})^{\rho_{n}}
|x_{11}x_{22}|^{\rho_{1,2}} \ldots |x_{1 n-m}\ldots x_{m
n}|^{\rho_{n-m,\ldots,n}}
$$
(with the condition $\lambda_s=\sum \rho_{i_1,\ldots,i_s}=l_s-l_{s+1}$,
$\lambda_k=l_k$) are not linearly independent.

Let us write down a $k$-connex $f$ of type $(l_1,\ldots,l_k)$
in the canonical form such that the size of its minor factors
decrease and the minors of the same size appear in the
lexicographical order,
assuming $1>2>\ldots>n$.
Every $k$-connex $f=f_1 f_2\ldots f_k$ ordered in such a way
(for $f_i= |x_{1 j_{1}^{i}}\ldots x_{l_i j_{l_i}^{i}}|$)
corresponds to a tableau
whose $i$-th column is formed by upright standing numbers
$j_{1}^{i},\ldots,j_{l_i}^{i}$. The tableau obtained has $k$
columns of corresponding lengths $l_1,\ldots,l_k$.
The entries of the tableau are numbers $1,\ldots,n$.

If a $k$-connex corresponds to a standard tableau then it is
called a standard $k$-connex.

\medskip

\textbf{Theorem 2 \cite{Hodge}.} {\it Standard $k$-connexes of type
$(l_1,\ldots,l_k)$ form a base of the vector space spanned by $k$-connexes
of this type. }

\begin{center}
\large 1.4 Chains
\end{center}

Let us put in order the basic vectors of $\mathbb{R}^n$ as
$e_1>e_2>\ldots >e_n$ and extend the order to basic vectors of
$\wedge^k \mathbb{R}^n$ lexicographically. Choose the canonical base
$\{e_{\alpha_1}\vee e_{\alpha_2}\vee \ldots \vee e_{\alpha_m} :
\alpha_1 \geq \alpha_2 \geq \ldots \geq \alpha_m\}$ in the space
$V(m,n,k)$. If a basic vector $v$ is in the form $v=v_1\vee\ldots
\vee v_m$ with lexicographically ordered components, then we call
such form {\it canonical}. If
$v=(v_{k_1})^{i_1}\vee\ldots\vee(v_{k_q})^{i_q}$, $v_{k_1}<v_{k_2}
<\ldots<v_{k_q}$, $i_1+\ldots +i_q=m$,
then let us denote
\begin{equation}\label{eq:CoeffL}
\lambda(v) = \frac{m!}{i_1!\ldots i_q!}.
\end{equation}
Introduce the scalar product $(..,..)$ on $V(m,n,k)$ supposing
$e_{\alpha_1}\vee e_{\alpha_2}\vee\ldots\vee e_{\alpha_m}$
is an orthonormal base.
Let us also denote a basic vector $e_{\alpha_1}\vee e_{\alpha_2}\vee \ldots \vee
e_{\alpha_m}$ by $(\alpha_1,\alpha_2,\ldots,\alpha_m)$.

It is well-known \cite{Kos} that
$$
\dim V(m,n,k)=\binom{m+t-1}{t-1}=\binom{m+t-1}{m},\
\mbox{where}\,\,t=\binom{n}{k}.
$$

For any vector $v$ from the canonical base of $V(m,n,k)$
define a set \linebreak
$(l_1,\ldots,l_n)_v$, where $l_i$ stands for how many times
a number $i$ occurs in the record of~$v$.

Let $(k_i)$, $i=1,\ldots,n$, be a set of numbers
$\mathbb{N}\cup\{0\}$, such that
\begin{eqnarray}  \label{Cond}
\sum_{i}^{n}k_i=mk,\quad k_i\leq m,\quad |\{i \,|\, k_i>0\}|\geq k.
\end{eqnarray}
A {\it chain} of type $(k_i)$ is a set of vectors $v$ from the canonical base
of $V(m,n,k)$ such that $(l_i)_v=(k_i)$.
It is clear that all chains form a partition of the standard base of $V(m,n,k)$.
Further, we will also mean by a chain a set satisfying the conditions (\ref{Cond}).

A vector $v$ of the canonical base of $V(m,n,k)$ is called {\it invariant},
if its corresponding chain contains only the vector $v$, i.e.,
$\{u\,|\, (l_i)_u=(l_i)_v\}=\{v\}$. Otherwise, it is called {\it non-invariant}.

Let us denote by $L(\sigma)=L(u_i\,:\, u_i \in\sigma)$
the linear envelope of a chain $\sigma=(k_i)$, and
let $Pr(v,\sigma)$ stands for the projection of a  vector $v\in V(m,n,k)$
onto the subspace $L(\sigma)$.
\medskip

\textbf{Lemma 1\cite{Gub2}.} \\
1)$\,\, V_0\,=\,L((x_1\wedge y_2\wedge \ldots \wedge y_k)\vee
\ldots \vee(x_m\wedge y_2\wedge \ldots \wedge y_k) : x_i, y_j \in
\mathbb{R}^{n})$,
\\
2) {\it $Pr(v,\sigma)\in V_0$ for every $v\in V_0$ and for every
chain $\sigma$.}

\medskip

\textbf{Theorem 3\cite{Gub2}.} {\it A vector from the
canonical base of
$S^{m}(\wedge^{k}\mathbb{R}^{n})$ belongs to $V_0$
if and only if it is invariant}.

\medskip

These results imply that a base and dimension of $V_0$
can be found as follows.
Let $\sigma$ be a chain
(hereinafter, by a chain we mean a chain consisting
of non-invariant vectors). If we can find a base
of all intersections $V_0\cap L(\sigma)$ (for all chains $\sigma =(k_i)$)
then the union of all these bases and all invariant vectors in $V(m,n,k)$
is a base of the entire space~$V_0$.

Let $\sigma=\{v_1,\ldots,v_t\}$ be an arbitrary chain,
$t=t(\sigma)=|\sigma|$, and let the vectors $v_i$ be ordered by decrease.
Let $w=x_1\wedge\ldots\wedge x_k,\>x_i=\sum_{j=1}^{n}\alpha_{ij}e_j$.
Consider
\begin{equation}
\label{Sistem} \left\{
\begin{array}{ccc}
A_1(\alpha_{ij}) & = & (w^m,v_1), \\
\ldots &  & \ldots \\
A_t(\alpha_{ij}) & = & (w^m,v_t),
\end{array} \right.
\end{equation}
where $A_s$ are homogenous polynomials of the power $mk$ in
variables $\alpha_{ij}$. By the definition, a polynomial $A_s$ can be
represented in the form of a product of $m$ minors of the size $k$
of the matrix $(\alpha_{ij})$, and the choice of minors is determined
by the  vectors $v_s$.

Since $w$ is decomposable, equations (\ref{Plue})
impose restrictions on matrix minors $(\alpha_{ij})$,
so $A_s(\alpha_{ij})$ also satisfy some relations.
Given a chain $\sigma $, denote by $r(\sigma)$
the rank of the system
$A_{1},\ldots, A_{t}$ (i.e., the maximal number of linearly
independent polynomials in this system).

\medskip

\textbf{Theorem 4\cite{Gub2}.} {\it For every chain $\sigma$
we have $\dim(V_0\cap L(\sigma))=r(\sigma)$.}

\medskip

\textbf{Remark 1.} In general case it is not hard to show an equivalence
of the problems of finding of base and dimension of the space $V_0$ and
analogous problem for the space of $m$-products of maximal matrices minors
$k\times n$ from variables $(x_{ij})$.

\medskip
In the case $m=2$, the following useful result was obtained in \cite{Gub2}.

\textbf{Lemma 2\cite{Gub2}.} {\it For the space $V(2,2k,k)$ and the chain
$\delta=(1,\ldots,1)$, a vector
$$
u=\sum_{\{i_1,\ldots,i_{2k}\}=\{1,\ldots 2k\}}q_{(i_1,\ldots i_k,
i_{k+1}\ldots i_{2k})} v_{(i_1,\ldots i_k,i_{k+1}\ldots i_{2k})}
$$
belongs to $V_0$ if and only if
\begin{eqnarray}\label{Soot3}
\sum_{s=1}^{k+1}(-1)^{s}q_{(i_1\ldots \dve{i_{s}}\ldots
i_{k+1},i_s i_{k+2}\ldots i_{2k})}=0
\end{eqnarray}
for arbitrary sets $1\leq i_1<\ldots<i_{k+1}\leq 2k$,
$1\leq i_{k+2}<\ldots<i_{2k}\leq 2k$, $i_c\not =i_d (c\neq d)$.}

\medskip
Let us denote (\ref{Soot3}) by $[i_{k+2}i_{k+3}\ldots i_{2k}]$, for short.

\medskip
\begin{center}
{\bf \large \S2 Base and dimension of $V_0$}
\end{center}

Suppose a basic vector $v=v_1\vee \ldots \vee v_m $ has ordered
on decrease components $v_i$. For this vector, consider the following
tabular representation:
$$
A(v) =
\begin{array}{|c|c|c|}
\hline a_{11} & \ldots & a_{1m} \\ \hline
 \ldots & \ddots & \ldots \\ \hline
 a_{k1} & \ldots & a_{km} \\ \hline
\end{array}
$$
for $v_j = e_{a_{1j}}\wedge \ldots \wedge e_{a_{kj}}$.

Clearly, the entries in each column strongly decrease.
If, in addition,
the entries in the rows weakly increase, i.e., $A(v)$
is standard then we name such vector standard.

Let us consider the space $U_0$ generated by Young diagrams
$(a_{ij})$ of the size $k\times m$ with the following relations:

$(S1)$ it is possible to interchange two entries in any column
with changing the sign of the table;

$(S2)$ it is possible to interchange every two rows;

$(S3)$ it is possible for any two chosen columns $j_1$, $j_2$ and
any row number $i_0$ to replace the primary table with the following
linear combination:
$$
A = \frac {1} {\lambda (A)} \sum _ {t=1} ^ {k} (-1) ^t\lambda (B_t) B _ {t},
$$
where $\lambda(C)=\lambda(v)$, see \eqref{eq:CoeffL},
 $v$ is the vector corresponding to table $C$, and
$$
(b_{t})_{i_0 j_1} = a_{t j_2},                      \quad
(b_{t})_{t j_2} = a_{i_0 j_1},                      \quad
(b_{t})_{ij} = a_{ij}, \mbox{otherwise}.
$$

Given a chain $\sigma$, consider a subspace $U(\sigma)\subseteq U_0$
generated by the images of all diagrams $A(v)$, $v\in\sigma$.
Relations $(S1)$, $(S2)$ provide correct one-to-one correspondence
between basic vectors and tables. Relations $(S3)$
correspond to Pl\"ucker equations.

Let us consider a ring $\mathcal R$ formally generated by minors
of the size $k$ from (\ref{Minors}). Pl\"ucker equations in a
minor form are homogeneous and quadratic, so they generate an ideal
$I$ in the ring $\mathcal R$.

Consider the $m$-th homogeneous component of the quotient ring  $\mathcal R/I$.
Images of two different $m$-products of minors $z_1$ and $z_2$ coincide
if and only if their difference is a linear combination
$f_1 P_1+\ldots+f_s P_s$, where $f_i$ are homogeneous polynomials on minors
of degree $m-2$, and $P_i$ are Pl\"ucker equations. Thereby, a problem of finding
a base and dimension of $U_0$ is equivalent to a similar problem for the space
of $m$-products of the maximum minors of a ($k\times n$)-matrix
$(x_{ij}) $.
By remark 1 it is equivalent to the initial problem on the space $V_0$.
Moreover, for every chain $\sigma$
 base and dimension search for spaces $V_0\cap L(\sigma)$
and $U(\sigma)$ are equivalent problems.

The reasoning stated above proves the theorem~7 \cite {Gub2} in
the inverse direction.

\medskip

A chain $\sigma$ is called {\em full} if
each index $i=1,\dots, n$ occurs no more than once
in the  record of all its vectors.
Then in the tabular representation
of every vector from $\sigma$ all entries are pairwise distinct.
The chain $\sigma=(1,1,\ldots,1,0,\ldots,0)$ with $mk$ units
is an example of a full chain.

Consider for any standard vector $v_s\in\sigma$ a vector
$$
t_s = Pr\big(\sigma, ((a_{11}^s + \ldots +a_{1m}^s)\wedge
(a_{21}^s + \ldots +a_{2m}^s)\wedge \ldots \wedge (a_{k1}^s +
\ldots +a_{km}^s))^m\big).
$$
By distributing the brackets we receive $t_s=\alpha v_s+w$, where $\alpha\not=0$,
$v_s>w$ in the lexicographic order. Thereby, vectors $t_s$ are linearly independent, and
$\dim U(\sigma)$ is greater or equals to the quantity of standard Young diagrams.

\medskip

Let us prove by an induction on $m$ that any non-standard tableau
corresponding to a vector from a full chain is linearly expressed
through standard ones that are higher than the primary tableau.

At $m=2$ it has been proved that  $\dim (V_0\cap L(\sigma))$
equals to the quantity of standard Young diagrams \cite{Gub2}.
Hence, any non-standard tableau corresponding to a vector from a full
chain is linearly expressed through standard ones.
We will show by induction on $k$,
that in this linear combination all standard vectors will be higher
than the vector corresponding to the primary tableau.
Moreover, in the rewriting process that leads to the desired linear combination
all relations of type $(S3)$ are applied with a choice of a number $i_0$ from the first column.

At $k=2$ and $k=3$ it can be
checked directly.
Let $k>3$, and let $u=(1 i_1\ldots i_{k-1}, j_1 \ldots j_k)$
be the vector corresponding to the primary table. If $i_1>3$ then by means
of relations $[i_1\ldots i_{k-1}]$ the vector $u$ can be expressed through
vectors of the form $(1 2 \ldots, \ldots)$ and $(1 3 \ldots, \ldots)$.

By the induction assumption for vector $(3 l_1\ldots l_{k-2},t_1 \ldots t_{k-1})$,
a vector of the form $(1 3 l_1 \ldots l_{k-2}, 2 t_1 \ldots t_{k-1})$
can be expressed as
a linear combination of standard ones and vectors of the form $(1 2 \ldots, \ldots)$.

Let us consider a vector $u=(1 2 i_1\ldots i_{k-2}, j_1 \ldots j_k)$.
Apply the induction assumption for vector $(2 i_1\ldots i_{k-2},j_2\ldots j_k)$:
Relations $(S3)$ allow to conclude that all vectors are higher than $u$.
Thereby, the statement is proved for $m=2$.

Assume the statement is proved for any $l<m$. One may order the
columns of the tableau in the decreasing order.
We will study the tableau formed
by the first $m-1$ columns. If it is not standard then we can linearly
express it through higher standard ones by the induction assumption.
If, after that, all tableaux have standard form then the statement is proved.
Otherwise, we express the tableau formed by its two last columns
through the higher standard ones. After that, we check
whether each of these whole tableaux belongs to the set of standard ones.
Note that last column becomes lower.

Further, considering and linearly expressing the tableaux formed by
first $m-1$ or last two columns we will stop at a final step.

\medskip

Consider a chain $\sigma = (k_i)$. Introduce a linear mapping
$\phi:\mathbb{R}^{mk} \rightarrow \mathbb{R}^n$ as follows
$$
\gathered \phi(e_s)=e_1, \ s=1,\ldots,k_1; \quad
\phi(e_s)=e_2, \ s=k_1+1,\ldots,k_1+k_2; \\
\ldots \\
\phi(e_s)=e_n, \ s=k_1+\ldots +k_{n-1}+1,\ldots,k_1+\ldots +k_n.
\endgathered
$$

The mapping $\phi$ can be extended to the base of $V$ and
thus it can be defined on tableaux. From the definitions of $\phi$ and
 $V_0$  it follows

\medskip

\textbf{Lemma 3.} {\it For every chain $\sigma$ and full chain $\sigma_0$
it is true that

1) $\phi(L(\sigma_0)\cap V_0) = L(\sigma)\cap V_0$;

2) $\phi(U(\sigma_0)) = U(\sigma)$. }
\medskip

Thereby,  $\phi$ maps standard tableaux
of $\sigma_0$ to standard tableaux of
$\sigma$, and the set of their images is a base of the space $L(\sigma)\cap V_0$.
We receive that the set of standard tableaux forms a base of $U(\sigma)$,
their total number can be found by the formula (\ref{StanleyF}):
\begin{multline*}
N(m,n,k) = \left( \frac{n!}{(n-k)!}
\frac{(n+1)!}{(n-k+1)!} \ldots \frac{(n+m-1)!}{(n+m-k-1)!}\right)\\
\times \left( \frac{1!}{k!} \frac{2!}{(k+1)!} \ldots
\frac{(m-1)!}{(k+m-1)!}\right) = \frac{\prod_{i=1}^m
\binom{n+i-1}{k}}{\prod_{j=1}^{m-1} \binom{k+j}{k}} =\left(
\prod_{i=1}^{m}\frac{\binom{n+i-1}{k}}{\binom{k+i}{k}} \right)
\binom{m+k}{k}.
\end{multline*}

For each standard vector $v$ with tabular representation
$A(v)$ we construct the following vector:
\begin{eqnarray}\label{Base}
\tilde{v} = \big((a_{11}+ \ldots +a_{1m})\wedge (a_{21}+\ldots
+a_{2m}) \wedge \ldots \wedge (a_{k1} + \ldots + a_{km})\big)^m.
\end{eqnarray}

Let us compile the results received in the following theorem
completely answering a question from \cite{Shar1}:
\medskip

\textbf{Theorem 5.} {\it Vectors built by (\ref{Base})
form a base of $V_0$, and
\begin{eqnarray}\label{Dim}
\dim V_0(m,n,k) = \left(
\prod_{i=1}^{m}\frac{\binom{n+i-1}{k}}{\binom{k+i}{k}} \right)
\binom{m+k}{k}.
\end{eqnarray}}

\begin{center}
{\bf\large\S3
Embedding the Grassmann variety into $S^{2}(\wedge^{k}\mathbb{R}^{n})$}
\end{center}

Consider a mapping $\phi: \wedge^k \mathbb{R}^n \mapsto V(2,n,k)$
defined by the rule $\phi(v)=(\sgn \,v) v^2$, where $\sgn \, v=\sgn
(v,e_{\alpha})$, $\alpha=\max \{\beta : (v,e_{\beta})\not =0\}$.
Denote an image $\phi(D)$ from the set $D$ of all nonzero decomposable
vectors of $\wedge^{k}\mathbb{R}^{n}$ as $D_0$.

Let us consider the following relations:
\begin{equation} \label{r1}
\begin{gathered}[]
 [j_1 j_2 i_1\ldots i_{s-3}], \ i_s\in \{j_3,\ldots,j_{2s}\},
\quad
 [j_1 j_3 i_1\ldots i_{s-3}], \ i_s\in \{j_4,\ldots,j_{2s}\},
\\
 [j_2 i_1\ldots i_{s-2}],\ i_s\in \{j_4,\ldots,j_{2s}\}.
\end{gathered}
\end{equation}
for any chain $\sigma =\sigma_v $, $v=(i_1\ldots
i_{k-s}j_{1}\ldots j_{s}, i_1\ldots i_{k-s} j_{s+1}\ldots
j_{2s})$;
\begin{eqnarray}
q_{\alpha,\beta}^2-4q_{\alpha,\alpha} q_{\beta,\beta}=0
\label{r2},
\\
q_{\alpha,\beta} q_{\alpha,\gamma}-2q_{\alpha,\alpha}
q_{\beta,\gamma}=0 \label{r3},
\end{eqnarray}
where $\alpha,\beta,\gamma$ are arbitrary pairwise distinct elements of $\Pi$.

\medskip
\textbf{Theorem 6.} Nonzero vector $u=\sum_{\alpha \geq
\beta}q_{\alpha,\beta}e_{\alpha,\beta}$ belongs to $D_0$ if and only if
the relations (\ref{r1})--(\ref{r3}) for coefficients $q_{\alpha,\beta}$
are satisfied.

\medskip

\textsc{Proof}.
Assume $0\ne u=\sum_{\alpha\geq\beta}q_{\alpha,\beta}e_{\alpha,\beta} \in D_0$.
Then Lemma 2 implies that
equations (\ref{r1}) hold.
Without loss of generality we may assume $u=v\vee v$ for some
$v=\sum_{\alpha\in \Pi}v_{\alpha}e_{\alpha}$. Then
$$
u=\bigg(\sum_{\alpha} v_{\alpha}e_{\alpha}\bigg)\vee
\bigg(\sum_{\beta}v_{\beta}e_{\beta}\bigg)=\sum_{\alpha}v_{\alpha}^2
e_{\alpha,\alpha}+2\sum_{\alpha >
\beta}v_{\alpha}v_{\beta}e_{\alpha,\beta},
$$
hence, relations (\ref{r2})--(\ref{r3}) hold.

Conversely, assume the relations (\ref{r1})--(\ref{r3}) hold for
coordinates of a non\-zero vector $u$. Let
$\alpha=\max\{\beta:q_{\beta,\beta}\not =0\}$. Without loss of
generality,  suppose that $q_{\beta,\beta}>0$. It is enough to show
that $u=v\vee v$. As $u\in V_0$ by (\ref{r1}), decomposability of
$v$ follows from Theorem 8 in \cite{Gub2}.

Construct a bijective mapping $o: \Pi \mapsto
\{1,2,\ldots,\binom{n}{k}\}$ such that
 $o(\alpha)>o(\beta)$ if and only if
$\alpha>\beta$.
Compose a matrix $A\in M_{\binom{n}{k}}(\mathbb{R})$ in the following way:
$$
a_{ij}= \left\{
\begin{array}{ll}
\nonumber q_{o^{-1}(i) o^{-1}(j)}, & i=j,\\
\nonumber \frac{1}{2} q_{o^{-1}(i) o^{-1}(j)}, & i\not=j.\\
\end{array} \right.
$$

The Lagrange method of reduction of a quadratic form $P(x)=\sum\limits_{i,j}a_{ij}x_{i}x_{j}$
to the canonical form implies that $u=\widetilde{e}_{1}^{2}$ in some base $\widetilde{e}_i$
if and only if $\rank\,A=1$.
Equivalently, all minors of size 2 of the matrix $Q$
are zero, i.e., relations (\ref{r2})--(\ref{r3}) hold as well as
relations
\begin{eqnarray}
q_{\alpha,\beta} q_{\gamma,\delta}-q_{\alpha,\delta}
q_{\gamma,\beta}=0 \label{r4},
\end{eqnarray}
where $\alpha,\beta,\gamma,\delta$ are arbitrary pairwise
distinct elements of $\Pi$.

And relations (\ref{r4}) follow from previous ones because of
\begin{multline}
(q_{\alpha,\beta} q_{\gamma,\delta}-q_{\alpha,\delta}
q_{\gamma,\beta})^2 = q_{\alpha,\beta}^2
q_{\gamma,\delta}^2+q_{\alpha,\delta}
q_{\gamma,\beta}^2-2q_{\alpha,\beta}q_{\gamma,\delta}q_{\alpha,\delta}
q_{\gamma,\beta} \\
= 8q_{\alpha,\alpha}
q_{\gamma,\gamma}(4q_{\beta,\beta}q_{\gamma,\gamma}-q_{\beta,\gamma}^2)=0.
\end{multline}

{\bf Corollary.}
It is possible to introduce the Grassmann variety $G_{n,k}$ in
coordinates of the space $S^{2}(\wedge^{k}\mathbb{R}^{n})$ as
the algebraic variety satisfying relations
(\ref{r1})--(\ref{r3}). Relations (\ref{r1}) are linear,
relations (\ref{r2}) and (\ref {r3}) are conic quadrics \cite{Vin}
of ranks 3 and 4, respectively.

\medskip

{\bf Remark 2.}
The number of relations (\ref{r1}) equals $\codim V_0(2,n,k)$,
of form (\ref{r2})--- $A=\binom{n}{k}\left(\binom{n}{k}-1\right)/2$, and of form
(\ref{r3})--- $\binom{n}{k}A$.

\begin{center}
{\bf\large
Acknowledgements}
\end{center}

The work is supported by the Federal Target Grant ``Scientific and
educational staff of innovation Russia'' for 2009--2013 (contracts
02.740.11.0429, and 14.740.11.0346).

The author expresses to V.A.Sharafutdinov gratitude for statement the problem,
P.S.Kolesnikov for a management of the given work.

\end{document}